\documentclass[10pt]{amsart}
\usepackage{amsmath, latexsym, 
}
\usepackage{amssymb,amsthm}
\usepackage{amscd,graphics}

\newtheorem{proposition}{Proposition}[section]
\newtheorem{lemma}[proposition]{Lemma}

\newtheorem{theorem}[proposition]{Theorem}

\newtheorem{corollary}[proposition]{Corollary}

\theoremstyle{remark}

\newtheorem{remark}[proposition]{Remark}

\theoremstyle{definition}

\def\real{\mathbb{R}}

\def\complex{\mathbb{C}}
\def\supp{\mathrm{supp}}
\def\var{\mathrm{var}}
\def\sgn{\mathrm{sgn}}

\def\id{\mathrm{id}}
\def\Lip{\mathrm{Lip}}

\def\BB{\mathcal{B}}
\def\DD{\mathcal{D}}

\def\LL{\mathcal{L}}
\def\MM{\mathcal{M}}

\def \PP{\mathcal {P}}

\def\UU{\mathcal{U}}
\def\VV{\mathcal{V}}

\begin{document}
\title[Alternative proofs of linear response]{Alternative proofs of linear response  for piecewise expanding unimodal maps}
\author{Viviane Baladi and Daniel Smania} 
\address{D.M.A., UMR 8553, \'Ecole Normale Sup\'erieure,  75005 Paris, France}
\email{viviane.baladi@ens.fr}

\address{
Departamento de Matem\'atica,
ICMC-USP, Caixa Postal 668,  S\~ao Carlos-SP,
CEP 13560-970
S\~ao Carlos-SP, Brazil}
\email{smania@icmc.usp.br}
\date{\today} 
\begin{abstract}
We give two new  proofs that the SRB measure $t \mapsto \mu_t$ of
a $C^2$ path $f_t$ of unimodal piecewise expanding $C^3$ maps is differentiable at $0$ 
 if $f_t$ is tangent to the topological class of $f_0$.
The arguments are more conceptual than the one in \cite{BS1}, but
require proving H\"older continuity of the infinitesimal conjugacy $\alpha$ 
(a new result, of independent interest) and using spaces of
bounded $p$-variation. 
The first new proof  gives differentiability
of higher order of $\int \psi \, d\mu_t$ if $f_t$ is smooth enough and stays in the
topological class of $f_0$ and if $\psi$ smooth enough
(a new result). In addition, this proof does not require any
information on the decomposition of the SRB measure into 
regular and singular terms, making it potentially amenable to
extensions to higher dimensions. The second new proof
allows us to recover the linear response formula (i.e., the formula
for the derivative at $0$) obtained in \cite{BS1}, and gives additional
information on this formula.
\end{abstract}
\thanks{V.B. is partially supported by ANR-05-JCJC-0107-01.
D.S. is partially supported by CNPq 470957/2006-9
and 310964/2006-7, FAPESP 2003/03107-9. 
D.S. thanks the DMA of Ecole Normale Sup\'erieure for
hospitality during a visit
where a crucial part of this work was done.}
\maketitle


\section{Introduction}

Many chaotic dynamical systems $f:M\to M$
on a Riemannian manifold $M$ admit an SRB measure $\mu$
(see e.g. \cite{Yo}) which  describes the statistical properties
of a ``large" set of initial conditions in the sense of Lebesgue
measure. (In dimension one, an SRB measure is simply an 
absolutely continous ergodic invariant
probability measure $\mu_t=\rho_t\, dx$ with a positive Lyapunov exponent.) It is of interest
(in particular in view of applications to statistical mechanics,
see e.g. \cite{Rua, Rub}) to study the smoothness of $t\mapsto \mu_t$,
when $f_t$ is a smooth family of dynamical systems, each having an SRB measure
$\mu_t$. If $t \mapsto \mu_t$ is
differentiable, one says that {\it linear response} holds.
Ruelle \cite{Rua} obtained not only differentiability, but also
a formula for the derivative (the {\it linear response formula}), 
in the case of smooth uniformly hyperbolic
dynamical systems.

In \cite{BS1}, we proved that the SRB measure
$t\mapsto \mu_t$ of a $C^2$ family of piecewise $C^3$ and
piecewise expanding unimodal maps
$f_t$, with $f_0$ mixing (see \S \ref{3.3} for formal definitions), is differentiable at $t=0$ (as a Radon measure) if and only
if $f_t$ is tangent to the topological class of $f_0$ at $t=0$.
(Keller \cite{Ke0} proved a long time ago that $\rho_t$ has a $|t ||  \ln |t||$ modulus
of continuity, as an element of $L^1(dx)$, and examples in \cite{Ba, BS1} show that 
this can be attained for non tangential families.)
We also obtained in \cite{BS1} a linear response formula analogous to the one in
\cite{Rua} (we used the resummation
introduced in \cite {Ba}).

More recently, differentiability of the SRB measure (in the weak sense,
that is, as an appropriate distribution)
was obtained \cite{RuStr}, \cite{BS3} for smooth families $f_t$
of analytic and nonuniformly expanding unimodal maps which stay in the topological
class of $f_0$.  The cases of families of
smooth nonuniformly expanding
interval maps only tangent  to the topological class
(where Whitney differentiability is expected on suitable subsets
of parameters),
as well as higher-dimensional dynamical systems such as piecewise
expanding/hyperbolic maps or H\'enon-like maps, are still open, and
much more difficult, see \cite{BN} for a discussion.
In particular, the arguments in \cite{BS1} and \cite{RuStr} used detailed
information about the structure of the SRB measure, decomposing it into a regular and a
singular term. This type of information may be far less accessible
in higher dimensions.

In this article, we give two new proofs of the fact 
(Theorem 5.1 in \cite{BS1}) that
the SRB measure of a $C^2$ family of piecewise $C^3$ and
piecewise expanding unimodal maps
$f_t$, with $f_0$ mixing,  is differentiable at $t=0$ if 
$f_t$ is tangent to the topological class of $f_0$ at $t=0$.

Section~\ref{firstt} contains our first new proof (see Corollary~\ref{newproof}), more
precisely, we obtain differentiability of $t \mapsto \int \psi \, d\mu_t$
for $\psi \in C^{1+\Lip}$ if $f_t$ is  a $C^2$ family of
piecewise expanding $C^3$ unimodal maps tangent to the topological class
of a mixing map $f_0$.
The argument is based on thermodynamic formalism, using
potentials $(s,t)\mapsto s (\psi \circ h_t) - \log |f'_t\circ h_t|$
(where $h_t$ conjugates $\tilde f_t$ with $f_0$ and $|\tilde f_t - f_t |= O(t^2)$)
and does not require any knowledge about
the structure of $\mu_t$.  It may therefore be useful in more difficult situations
(such as H\'enon maps, see \cite{BN}).
It requires proving H\"older differentiability of the
{\it infinitesimal conjugacy} $\alpha$,
a new result (Proposition ~\ref{glu}), of independent interest.
Also, this new proof gives  that $t \mapsto \int \psi d\mu_t$ is
a $C^j$ function, if $\psi \in C^{j+\Lip}$ and
$f_t$ is a $C^{j+1}$ family of piecewise expanding $C^{j+2}$ maps
in the topological class of $f_0$,
for any $j\ge 1$ (this is a new result, Theorem~\ref{pressuretrick}).
Note also that we do not require the assumption from \cite{BS1} that there
is a function $X$ so that $\partial f_t|_{t=0}= X\circ f_0$.

The first new proof requires $\psi \in C^{1+\Lip}$
(instead of $\psi\in C^0$ as in \cite{BS1}) and  does
not furnish the linear response formula. 
Section ~\ref{secc} contains our second new proof
(Theorem~\ref{theformula}),
which  uses spectral perturbation theory for transfer operators
associated to the dynamics $f_0$ and the weight $1/|f'_t\circ h_t|$.
This other proof
gives differentiability of $\int \psi\, d\mu_t$ for $\psi \in C^0$ and, using
the assumption that $\partial f_t|_{t=0}= X\circ f_0$,
allows us to recover the linear response formula from \cite{BS1}.
(This second proof also uses the H\"older regularity of $\alpha$ from
Proposition~\ref{glu}.)
Note however that this second proof requires information on the
structure of $\mu_t$ from \cite[Prop. 3.3]{Ba}.

Putting together Theorem~\ref{pressuretrick} and ~\ref{theformula}
we get  the following additional result (Corollary~\ref{thebest}):
If $f_t$  is a $C^2$ family of piecewise
expanding $C^3$ unimodal maps
in the topological class of $f_0$, then 
$t \mapsto \mu_t$ is $C^{1}$ from a
neighbourhood of zero to Radon measures.

We emphasize that neither new
proof 
gives that the condition to be tangent to the topological class
is necessary, contrary to the argument in \cite{BS1} (see Theorem 7.1 there).
The proofs here are a bit shorter than the
one given in \cite{BS1}, although the present account 
requires some results from our previous papers
(such as \cite[Prop. 3.3]{Ba}
 \cite[Prop. 3.2, Theorem 2]{BS2} 
 \cite[Prop. 2.4, Lemma 2.6, Prop. 3.3]{BS1}).

\pagebreak

\section{Definitions and notations -- H\"older smoothness of the infinitesimal conjugacy}

\subsection{Formal definitions}
\label{3.3}

\medskip

Denote $I=[-1,1]$. 
For an integer $k \ge 1$, we define 
$\BB^{k}$  to be
the linear space   of  continuous functions 
$f\colon I \to\real$ such that $f$ is $C^{k}$ on the intervals 
$[-1,0] $ and  $[0,1]$.
Then $\BB^{k}$ is a Banach space for the norm
$ \max \{|f|_{C^{k}([-1,0])}, |f|_{C^{k}([0,1])}  \}$.
For an integer $k\ge 1$, we define the set  $\UU^k$  of
{\it piecewise expanding $C^ k$ unimodal maps}
to be the set  of   $f \in \BB^k$ such that 
\footnote{A prime denotes differentiation with respect to  $x\in I$,
a priori in the sense of distributions.}
$f(-1)=f(1)=-1$, $\inf_{x \ne 0}  |f'(x)| > 1$,
and  $f(0)\leq 1$
(so that $f(I)\subset I$). The point $c=0$ is called the {\it critical point} of $f$.

A piecewise expanding $C^{k}$ unimodal map $f$ is {\it good} if
either $c$ is  not periodic under $f$
or  $|(f^{q-1})'(f(c)) |\min \{ |f_+'(c)|,|f_-'(c) |\} > 2$, where $q\ge 2$ is the minimal
period of $c$; it is  {\it mixing} if $f$ is topologically mixing on
$[c_2, c_1]$, where $c_k=f^k(c)$.

For $1\le j \le k$, a {\it $C^ j$ family of piecewise expanding $C^k$ unimodal
maps}  is  is a $C^j$ map $f_t$ from $t\in (-\epsilon, \epsilon)$ to $\UU^k$
for some $\epsilon >0$.
(In this paper, $k\ge 1$ is an integer
and $j$ is either an integer or $j=k-1+\Lip$ for $k\ge 2$, 
the notations $\BB^{k+\Lip}$ and
$\UU^{k+\Lip}$ for integers $k\ge 1$  being self-explanatory.
See also Remark~\ref{rkhol}.)

\begin{remark}\label{famdef}
A $C^j$ family $f_t$ of piecewise expanding $C^k$ unimodal maps
is a $C^{j,k}$ perturbation of $f_0$ in the sense of \cite{BS1} if $j= k\ge 2$.
\end{remark}

\begin{remark}\label{rkhol} Considering  $\BB^{k+\beta}$
and $\UU^{k+\beta}$ for $k\ge 1$ integer and a H\"older exponent
$0<\beta<1$ will perhaps  allow to avoid the loss of regularity
from $C^{k+1}$ to $C^{k-1+\Lip}$ e.g. in \cite[Prop. 2.4]{BS1}
(this question was asked by J.-C. Yoccoz).
However, since the spectral result of Wong \cite{W} 
only holds on the space $BV_p$ of functions
of bounded $p$-variation if $1\le p <p_0$,
for some $p_0>1$ depending on the dynamics,
it may be necessary in this case to replace
$BV_p$ by  spaces of generalised $p$-variation, as introduced
by Keller \cite{Ke}. (See also Remark~\ref{whyBVp}.)
\end{remark}

Assume that $f_t$ is a $C^j$ family
of piecewise expanding $C^k$ unimodal maps for $k\ge j>1$. 
By classical results of Lasota--Yorke, 
each $f_t$
has a unique absolutely continuous invariant probability measure
$\mu_t=\rho_t\, dx$. This measure is ergodic and it is called the
{\it SRB measure} of $f_t$. If $f_t$ is mixing, then $\mu_t$ is mixing.
If $f_0$ is good and mixing, then  $f_t$ is mixing for all small
enough $t$ (see \cite{Ke0} and references therein).

We say that a  $C^ j$ family $f_t$ of piecewise expanding $C^ k$ unimodal
maps    is in the {\it topological class of $f_0$} if 
there exist homeomorphisms $h_t:I\to I$ such that
\begin{equation}\label{starr}
h_0=\id \mbox{ and } h_t\circ f_0=f_t\circ h_t\, , \quad \forall \, |t|<\epsilon\, .
\end{equation}
We say that a $C^j$ family $f_t$ of piecewise expanding $C^k$ unimodal
maps  ($k\ge j \ge 2$) is  {\it tangent to the topological class}
of $f_0$ if there exists a $C^{j}$ family  $\tilde f_t$
of piecewise expanding $C^k$ unimodal maps in the topological class of $f_0$
so that $\tilde f_t=f_0$ and
$\partial_t f_t|_{t=0}=\partial_t \tilde f_t|_{t=0}$.
(Note that there is a typographical mistake in \cite[p. 682, line 6]{BS1}, where
``$C^{2,2}$ perturbation" should be replaced by ``$C^{r_0,r}$ perturbation.")

We say   that a bounded function $v:I \to \real$ is {\it horizontal} for $f$, if
$v(-1)=v(1)=0$, and 
setting $M_f=q$ if $c$ is periodic of minimal period $q$,
and $M_f=+\infty$ otherwise,
\begin{equation}\label{Mf}
J(f,v)= \sum_{j = 0}^{M_f-1} \frac{v(c_j)}{(f^{j})'(c_1)}=0 \, .
\end{equation}

In \cite[Cor. 2.6]{BS1} we proved that if $f_t$ is a $C^2$ family
of piecewise expanding $C^2$ unimodal maps
tangent to the topological class of $f_0$, 
then $v=\partial f_t|_{t=0}$ is 
horizontal  for $f_0$.
By \cite[Theorem 2]{BS1},
if  $f_t$ is a $C^2$ family
of piecewise expanding $C^2$ unimodal maps with $f_0$ good and
 $v=\partial f_t|_{t=0}$  is $C^2$ and
horizontal  for $f_0$, then
 there exists a $C^{1+\Lip}$ family  $\tilde f_t$
of piecewise expanding $C^2$ unimodal maps in the topological class of $f_0$
so that $\tilde f_t=f_0$ and
$\partial_t f_t|_{t=0}=\partial_t \tilde f_t|_{t=0}$.

We proved in \cite[Lemma 2.2]{BS1} that if $v:I\to \real$ is 
bounded  then
the twisted cohomological equation (TCE) 
\begin{equation}\label{tceeq}
v(x)=\alpha(f(x))- f'(x) \alpha(x)\, , \quad \forall x \in I \, , \quad x\ne c\, ,
\end{equation}
admits a  unique bounded solution $\alpha$ satisfying 
 $\alpha(c)=0$. This solution is obtained
as follows: If $c$ is not in the forward orbit of $x$,  set $M(x)=\infty$ and
 otherwise let $M(x)$ be the smallest integer $j\ge 0$ satisfying $f^j(x)=c$,
then put
\begin{equation}\label{bdd}
\alpha(x)=
\sum_{i=0}^{M(x)-1} \frac{v(f^i(x))}{(f^i)'(x)} \, .
\end{equation}
The function $\alpha$
is called the {\it infinitesimal conjugacy.} 

If $u: I \to \real$ is H\"older we denote its H\"older norm by
$|u|_\beta$.
Slightly abusing notation, we shall sometimes write
$ \partial_t f_t$ for  $\partial_s f_s|_{s=t}$, and similarly for other
functions depending on $t$.

\subsection{H\"older smoothness of the infinitesimal conjugacy $\alpha$}
\label{holds}

A  new result  that we shall require throughout
(see  Lemmas~\ref{jensenlemma} and ~\ref{gla}) is:

\begin{proposition}[Smoothness of the infinitesimal conjugacy]\label{glu}
Let $f\in \UU^2$ be such that $c$ is not
periodic. For any  $\beta \in (0,1)$  
there exist $C_\beta>0$ and  $\VV_\beta$ 
a neighbourhood
of $f$ in $\UU^2$ so that, for any $g\in \VV_\beta$ and every
$\beta$-H\"older $v:I\to \real$, with 
$v(-1)=v(1)=0$  and  $J(g,v)=0$, 
the unique bounded  function  $\alpha$  
(\ref{bdd}) satisfying $\alpha(c)=0$
and 
$v(x)=\alpha(g(x))- g'(x) \alpha(x)$ for all $x\ne c$
is $\beta$-H\"older,
with  
$$|\alpha|_{\beta}
\le C_\beta |v|_{\beta} \, .
$$

If the critical point of
$f\in \UU^2$ is periodic, the statement holds up to taking 
(for appropriate $\xi(\beta)>0$)
$$\VV_\beta=
\{ g 
\mid
\|g-f\|_{\BB^2}
< \xi(\beta)\, ,
 \exists \mbox{ homeomorphism } h:I\to I \mbox{ s.t. } g\circ h= h \circ f\}\, .
$$
 \end{proposition}

In particular, if $f \in \UU^2$  and $J(f,v)=0$ for
some Lipschitz $v$ with $v(-1)=v(1)=0$, 
the function $\alpha$ solving (\ref{tceeq}) is $\beta$-H\"older for any 
$\beta <1$.

\begin{remark} J\'er\^ome Buzzi \cite{Bu} showed us a simple
proof that if $h$ is a homeomorphism
so that $h\circ f = g\circ h$, for two piecewise
expanding $C^1$ unimodal maps $f$ and $g$, then
$h$ is $\beta$-H\"older, 
for any $\beta < \log  (\inf |g'|/2)/\log (2\sup |f'|)$. 
This fact neither implies nor is implied by Proposition~\ref{glu}.
\end{remark}

\begin{proof} 

{\bf Step I.}
For any $\beta<1$, there exist 
a neighbourhood $\VV_\beta$  of $f$ in $\UU^2$, $\ell \ge 1$, and $\eta >0$
such that $\lambda=(\inf_{g \in \VV_\beta}
\inf_{x \ne c} |g'(x)|)^{-1}<1$, and, for any
$g \in \VV_\beta$,
letting $d_1 < d_2 < \dots <d_p$ be the critical points of $g^\ell$,
putting $d_0=-1$, $d_{p+1}=1$, and setting
\begin{equation}
\label{thetadef}
\theta= \max_{0\le i \le p}
\sup_{\substack{x,y \in (d_i, d_{i+1}) \\ |x-y|< \eta}}
\frac{ |(g^\ell)'(x)|^\beta } { |(g^\ell)'(y)|} \, ,
\end{equation}
we have $2 \theta <1$.

Put
$\Delta_g=\min_{0\le i \le p} \{ d_{i+1}-d_i  \}$. Then
$\inf_{g \in  \VV_\beta} \Delta_g>0$ if  the critical point of $f$ is
not periodic.
Otherwise we have $\inf_{g \in  \VV_\beta} \Delta_g>0$, up to
replacing $\VV_\beta$ by a $\BB^2$-neighbourhood of $f$ 
in its topological class. In particular, we can assume
that $\eta < \inf_{g \in  \VV_\beta} \Delta_g$.
From now on, we fix   $\VV_\beta$, $\ell\ge 1$, and $\eta >0$ as above.

{\bf Step II.}
We claim
it suffices  to show the lemma for $g \in \VV_\beta$  with a periodic critical point:  
Indeed, if $g$  has a nonperiodic critical point, then
we consider $g_t =g + tw$ with $g_t\in \UU^2$,
$w \in \BB^2$, $w(-1)=w(1)=0$, and $J(g,w)\ne 0$. By \cite[Corollary 4.1]{BS2},
there exists a sequence $t_n\to 0$ such that each $g_n=g_{t_n}$ has a periodic critical 
point. 
In particular, $g_n$ converges to $g$ in the $\UU^{2}$ topology.
Then,
 by \cite[Proposition 3.2]{BS2}
  we have
$\lim_{n \to \infty} J(g_n,v)=0$.  Let  $w_n$ be
a $\beta$-H\"older function, with 
$w_n(-1)=w_n(1)=0$ and $|w_n|_\beta \le 1$, such that $J(g_n,w_n)=1$. Set
$$v_n=   v - J(g_n,v)w_n\, .
$$ 
Then we have $J(g_n,v_n)=0$
and $\lim_{n\to \infty} |v_n -v|_\beta =0$.
If  the proposition
holds for maps in $\VV_\beta$ with a periodic turning point,  
the unique function $\alpha_n$
so that $\alpha_n(c)=0$
and $v_n(x)=\alpha_n(g_n(x))- g_n'(x) \alpha_n(x)$ for all $x\ne c$,
satisfies 
$|\alpha_n|_\beta \leq C_\beta |v_n|_\beta$.
We can choose a subsequence $\alpha_{n_i}$ 
converging in the sup norm to a function $\alpha$. If follows from the uniform convergence of $\alpha_{n_i}$ 
that $\alpha$ satisfies the TCE (\ref{tceeq}) 
for $g$ and $v$,
and  that $|\alpha|_\beta \leq C_\beta |v|_\beta$.

{\bf Step III.} 
We assume from now on that $g\in \VV_\beta$ has a periodic turning point. 
The proof will be via an ``infinitesimal pull-back" argument.

First, since $J(g,v)=0$, it  is easy to see that there exists a 
$\beta$-H\"older function $\alpha_0:I \to \real$
with $\alpha_0(-1)=\alpha_0(1)=\alpha_0(c)=0$,
$\alpha_0(g(c))=v(c)$, and
 $$
 v(x)=\alpha_0(g(x))
-g'(x) \alpha_0(x)\mbox{  for every $x\ne c$ in the (finite) forward
orbit of $c$.}
$$

Second, we define by induction 
 continuous functions  $\alpha_i: I\to \real$, for $i\ge 1$,
such that $\alpha_i(-1)=\alpha_i(1)=\alpha_i(c)=0$,
$\alpha_i(g(c))=v(g(c))$,  that
 \begin{equation}\label{Ai}
 v(x)=\alpha_i(g(x))-g'(x) \alpha_i(x) 
 \mbox{ for every $x\ne c$ in the (finite) forward
orbit of $c$,}
\end{equation}
and, in addition, 
\begin{equation}\label{Di}
v(x)=\alpha_{i-1}(g(x))- g'(x)\alpha_i(x) \, , \forall x \neq c\, .
 \end{equation}
Indeed, suppose we have defined $\alpha_i$, for $0\le i\leq n$. Set
$\alpha_{n+1}(c)=0$, and
$$\alpha_{n+1}(x) =  \frac{\alpha_n(g(x))-v(x)} {g'(x)}\, ,\quad x \neq c\, .
$$
 Clearly, $\alpha_{n+1}(-1)=\alpha_{n+1}(1)=0$,  and  (\ref{Di}) holds for  $i=n+1$.  
Thus, since $v(x)=\alpha_n(g(x))-g'(x) \alpha_n(x)$ for every $x\ne c$ in the forward
orbit of $c$,  we find $\alpha_n(x)=\alpha_{n+1}(x)$ for each $x\neq c$ in the
forward orbit of $c$. 
Since $\alpha_{n}(c)=\alpha_{n+1}(c)=0$, we conclude that 
(\ref{Ai}) holds for $i=n+1$,   
and 
$\alpha_{n+1} (g(x))=v(x)$.  Last, but not least,   $\alpha_{n+1}$ is continuous 
on $I$ because $\alpha_{n} (g(x))=v(x)$.

Thirdly,  if $x$ is not a critical point of $g^j$, we set $v_0(x)= 0$, and
$$v_j(x)= \sum_{i=0}^{j-1} (g^{j-1-i})'(g^{i+1}(x))v(g^i(x)) \, ,
\quad j\ge 1 \, .$$
Recalling the notation $\ell$, $\{d_i\}$, from Step I, it is easy to see that
\begin{equation}\label{tceell}
v_\ell(x) =\alpha_{j\ell}(g^\ell(x)) -(g^\ell)'(x)\alpha_{(j+1)\ell} (x)   \, ,
\quad \forall x \notin \{ d_1, \dots d_p\} \, , \forall j \ge 0 \, .
\end{equation}
For $\ell \ge 2$
the function $v_\ell$  may have jump discontinuities  at the critical
points $d_i$ of $g^\ell$, but it is $\beta$-H\"older in the connected components
of $I\setminus \{d_1, \ldots, d_p\}$.

Finally, we shall use the  iterated twisted cohomological
equation (\ref{tceell}) to show that  there exists $C_\beta <\infty$
so that, for all $g \in \VV_\beta$ with a periodic
turning point and all
$\beta$-H\"older $v$ with $J(g,v)=0$  
(and $v(-1)=v(1)=0$), there exists $j_0\ge 0$  so that
$
|\alpha_{j\ell}|_\beta \le C_\beta |v|_\beta$ for all $j \ge j_0$.
In view of this, for $j\ge 0$, set $K_j^0=\sup_I |\alpha_{ j\ell}|$, 
$L^0= \sup_I |v_\ell|$, 
$$
K_j^\beta= \sup_{x\ne y}
\frac{|\alpha_{j\ell} (x)- \alpha_{j\ell}(y)|}
{|x-y|^\beta} 
\, , \, 
\widehat K_j^\beta=
\max_{0\le i\le p} \sup_{\substack{ x\ne y,  |x-y|<\eta \\ x,y \in ( d_i,  d_{i+1})}}
\frac{|\alpha_{j\ell} (x)- \alpha_{j\ell}(y)|}
{|x-y|^\beta} \, ,
$$
and
$$
L^\beta=
\max_{0\le i\le p} \sup_{\substack{x\ne y \\ x,y \in (d_i, d_{i+1})}}
\frac{|v_{\ell} (x)- v_{\ell}(y)|}
{|x-y|^\beta} \, ,
\quad D=  \max_{0 \le i\le p} \sup_{x,y \in (d_i, d_{i+1})}
 \frac{ |(g^\ell)''(x)|}{ |(g^\ell)'(y)|^2}
\, .
$$
Clearly, $\max(L^0,L^\beta) \le \widetilde C_\beta |v|_\beta$ for
all $g$ and $v$ under consideration, and we have
\begin{equation}\label{clear}
\frac{|\alpha_{(j+1)\ell} (x)- \alpha_{(j+1)\ell}(y)|}
{|x-y|^\beta}\le 2  \eta^{-\beta} K^0_{j+1} \quad \mbox{ if } |x-y|\ge \eta
\,  , \, \forall j \ge 0 \, .
\end{equation} 
Therefore, recalling the definition of $\lambda$ and (\ref{thetadef})
from Step ~ I, it suffices to show that
\begin{equation}\label{cclaim}
K_{j+1}^0 \le \lambda^\ell( K_{j}^0 + L^0)\mbox{ and }
\widehat K_{j+1}^\beta \leq (L^0+ K^0_j)D  + \lambda^\ell
L^\beta +\theta K_{j}^\beta\, ,\, \forall j \ge 0\, .
\end{equation}
Indeed, continuity of $\alpha_{(j+1)\ell}$ together with
 (\ref{clear}) (recall
also that $\eta < \inf_g \Delta_g$,
so that if $|x-y|< \eta$ then $[x,y]$ contains at most one
point $d_i$) imply
$$ K_{j+1}^\beta\le \max(2  \eta^{-\beta} K^0_{j+1} ,2\widehat K_{j+1}^\beta)\, .
$$
The above bound together with (\ref{cclaim})  yield
$E_\beta<\infty$ so that, for all $g\in \VV_\beta$ with a periodic turning point,
and all
$\beta$-H\"older $v$ with $J(g,v)=0$,
$v(-1)=v(1)=0$, there exists
$j_0\ge 0$ so that
\begin{equation*}
K_{j+1}^\beta \leq  E_\beta |v|_\beta + 
2\theta K_{j}^\beta\, ,\, \forall j \ge j_0\, .
\end{equation*}
Since $2 \theta <1$,  we  conclude by  a geometric series, taking
larger $j_0$ if necessary.

It remains to show (\ref{cclaim}).
We concentrate on the second bound
(the first is easier and left to the reader).
Let  $x, y \in (d_i, d_{i+1})$ satisfy $|x-y|< \eta$. Then  (\ref{tceell}) implies
(since $g \in \UU^2$, 
the function $(g^\ell)'$ is $C^1$ in the intervals of monotonicity of $g^\ell$)
\begin{align*}
&|\alpha_{(j+1)\ell}(x)-\alpha_{(j+1)\ell}(y)| 
\le |v_\ell(x)+\alpha_{j\ell}(g^\ell(x)) 
|  \Big|  \frac{1}{(g^\ell)'(x)} -  \frac{1}{(g^\ell)'(y)}\Big| \\
&\qquad \qquad \qquad\qquad\qquad\qquad\quad+  
\frac{|v_\ell(x)-v_\ell(y)|+|\alpha_{j\ell}(g^\ell(x))-\alpha_{j\ell}(g^\ell(y))|  }{|(g^\ell)'(y)|} \\
&\qquad \leq 
(L^0+ K^0_j) \frac{\sup_{[x,y]} |(g^\ell)''|}{\inf_{[x,y]} | (g^\ell)'|^2}   |x -y|  
 +  
\frac{L^\beta+\sup_{[x,y]} |(g^\ell)'|^\beta K_j^\beta} {\inf_{[x,y]} |(g^\ell)'|}|x-y|^\beta      \\
&\qquad \leq 
(L^0+ K^0_j) D   |x -y|   +  
(\lambda^\ell L^\beta+\theta K_j^\beta) |x-y|^\beta  \, .
\end{align*}

{\bf Step IV.}
Defining
$\tilde{\alpha}_n=\frac{1}{n}  \sum_{j=0}^{n-1} \alpha_{j\ell}$,
we can choose a subsequence $\tilde{\alpha}_{n_i}$ 
converging uniformly on $I$ to a function $\tilde \alpha$ satisfying
$|\tilde \alpha|_\beta \leq {C}_\beta|v|_\beta$.
By (\ref{tceell}), 
\begin{equation}\label{q1}
 v_\ell(x)= \tilde \alpha(g^\ell(x))-(g^\ell)'(x) \tilde \alpha (x)\, , \quad
\forall x \notin \{ d_1, \dots d_p\}
 \, .
\end{equation}
Let $\alpha\colon I\to \real$ be the unique bounded solution vanishing at $c$
to the TCE (\ref{tceeq}) for $g$ and $v$, as in (\ref{bdd}). Then
\begin{equation}\label{q2}
v_\ell(x)=  \alpha(g^\ell(x))-(g^\ell)'(x) \alpha(x)  \, ,
\quad \forall x \notin \{ d_1, \dots d_p\}\, .
\end{equation}
Since $\tilde{\alpha}$ is continuous   (\ref{q1}) and (\ref{q2}),
imply  $\alpha= \tilde{\alpha}$. We proved
$|\alpha|_\beta \le C_\beta |v|_\beta$,
for all $g \in \VV_\beta$ 
with a periodic turning point, and thus the proposition.
  \end{proof}


\subsection{Banach spaces of bounded variation}

We shall consider the Banach space of functions of bounded variation
$$
BV
=\{ \varphi : \real \to \complex \mid
\var (\varphi) < \infty\, , \supp (\varphi) \subset I\} / \sim\, ,
$$ 
endowed with the norm $\| \varphi\|_{BV} =\inf _{\psi \sim \varphi} \var (\psi)$,
where $\var$ denotes  total
variation, and $\varphi_1 \sim \varphi_2$ if the
bounded functions $\varphi_1$, $\varphi_2$ 
differ on an at most countable set. 
In addition, for $1 \le p < \infty$ we shall  work with
the Banach space of functions of bounded $p$ variation
(used in interval dynamics by Wong \cite{W})
$$
BV_p
=\{ \varphi : \real \to \complex \mid
\var_p (\varphi) < \infty\, , \supp (\varphi) \subset I\} / \sim\, ,
$$
where
$$
\var_p (\varphi)=
\sup_{x_1 < x_2 < \ldots < x_n } \biggl (\sum_{i=1}^n  |\varphi(x_{i+1})- \varphi(x_i)|^p
\biggr )^{1/p}\, ,
$$
the supremum ranging over all ordered finite subsets of $\real$.
Note that $\var_1=\var$ and $BV=BV_1$.
Wong \cite{W} does not quotient by the equivalence relation
$\varphi_1 \sim \varphi_2$,
but his results
remain unchanged if we consider elements in
$BV_p$ modulo  $\sim$ (a function in $BV_p$ is continuous except on an
at most countable set, see also
\cite[Lemma 1.4.a, Lemma 2.7]{Ke} and \cite{Bru}).
Note that for each $p \ge 1$ there is $C \ge 1$
so that $|\varphi|_\infty \le C \|\varphi\|_{BV_p}$ for all $\varphi$,
and  if $\varphi $ is $1/p$-H\"older, then
$ \| \varphi \|_{BV_p} \le   | \varphi |_{1/p}$.
Also,
\begin{equation}\label{algebra}
\| \varphi_1 \varphi_2\|_{BV_p} \le 2 \| \varphi_1\|_{BV_p} \| \varphi_2\|_{BV_p} \, ,\, \, 
\forall p \ge 1 \, ,
\end{equation}
and $\|\varphi \circ h\|_{BV_p}= \|\varphi \|_{BV_p}$
for any homeomorphism $h:I\to I$ and all $p\ge 1$.

\begin{remark}\label{whyBVp}
The reason we consider  spaces $BV_p$ for $p\ne 1$
 is because we are
concerned with differentiability in the $t$-parameter and we will
have to deal with derivatives
$\partial_t (\psi\circ h_t) |_{t=0}= \psi'  \cdot \partial_t h_t|_{t=0}$
or $\partial_t (f'_t \circ h_t)|_{t=0}= f''_0 \partial_t h_t|_{t=0}+v'$,
where $v'=\partial_t f'_t|_{t=0}$ is $C^1$, but 
$\partial_t h_t|_{t=0}$
does not belong to $BV$ in general.  We shall see, however, that
 Proposition~ \ref{glu} implies that $\alpha=\partial_t h_t|_{t=0}$ lies in
$BV_p$ for   all $p >1$.
\end{remark}


\section{Weak differentiability of the SRB via the pressure}
\label{firstt}

The main result of this section (Theorem~\ref{pressuretrick})
says that  for any
$j\ge 1$,  if $f_t$ is a $C^{j+1}$ family of piecewise expanding $C^{j+2}$
unimodal maps  in  the topological class of $f_0$,
then $R(t)=\int \psi\, d\mu_t$ is $C^{j}$ if $\psi$ is
$C^{j+\Lip}$. 
Even if $j=1$, this is a new result (Theorem 5.1 in \cite{BS1} only
gives differentiability at $t=0$).
The  argument is based on the topological pressure of
the potential $(s,t)\mapsto - \log |f'_t \circ h_t| + s (\psi \circ h_t)$
for the map $f_0$.
It is simple, but  does not give the formula for $\partial_t R(t)|_{t=0}$
(or higher order derivatives).
Using the linear response
formula from \cite[Theorem 5.1]{BS1}
or  Theorem ~\ref{theformula} below, Theorem~\ref{pressuretrick}
will imply Corollary~\ref{thebest}.

\begin{theorem}\label{pressuretrick}
For any  integer $j\ge 1$, if  $f_t$ is a $C^{j+1}$ family
of piecewise expanding $C^{j+2}$  unimodal maps in the
topological class of a mixing map
$f_0$, then there is $\hat \epsilon>0$
so that for any  $C^{j+\Lip}$  function $\psi$
the map $R(t)=\int \psi \rho_t \, dx$ is $C^j$ in
$(-\hat \epsilon,\hat \epsilon)$.
\end{theorem}

As an immediate corollary of Theorem~\ref{pressuretrick} and Proposition~\ref{KLbd},
we recover the first claim of \cite[Theorem 5.1]{BS1} if $\psi$ is
$C^{1+\Lip}$ (we do not need the assumption
$\partial f_t|_{t=0}=X \circ f_0$ used in \cite[Theorem 5.1]{BS1}):

\begin{corollary}\label{newproof}
Assume that $f_t$ is a $C^2$ family
of piecewise expanding $C^3$  unimodal maps, where $f_0$ is a good mixing map.
If $f_t$ is tangent to the topological class of  $f_0$
then for any $C^{1+\Lip}$ function $\psi:I \to \complex$, the map
$R(t)=\int \psi \, d\mu_t $
is differentiable at $t=0$.
\end{corollary}

\begin{proof}[Proof of Theorem~\ref{pressuretrick}]
Fix $\psi \in C^{1+\Lip}$, recall the notation $h_t$ from (\ref{starr}), put
\begin{equation}\label{gst}
g_{s,t}(y)=  \frac{\exp( s \psi( h_t(y))}{|f'_t ( h_t(y))|}\, ,
\quad y \in I \setminus\{c\}\, ,
\end{equation}
and consider the transfer operator
\begin{equation}\label{xfer1}
\widetilde \LL_{s,t} \varphi (x)=\sum_{f(y)=x} g_{s,t}(y) \varphi(y)\, .
\end{equation}
Note that $\LL_0=\widetilde \LL_{0,0}$ is the usual transfer operator for $f_0$.
We let $\LL_0$ act  on $BV_p$  for any fixed  $p \in [1, p_0)$, where 
$p_0$, depending on $f_0$ through $\inf |f_0'|$ and $\sup |f_0'|$, is given by the
main Theorem of \cite{W}, which says that $\LL_0$ on $BV_p$ has spectral
radius $1$, essential spectral radius $<1$, and $1$ is the only eigenvalue
of modulus $1$. Furthermore, $1$ is a simple eigenvalue with
an eigenvector $\rho_0$ which is strictly positive on
$[c_2, c_1]$. The corresponding fixed vector $\nu_0$ of $\LL^*_0$ is just Lebesgue
measure $dx$. We normalise so that $\int \rho_0 d\nu_0=1$ and
$\nu_0(I)=1$. (Of course, $\mu_0=\rho_0 \, dx$.)
We shall view $\widetilde \LL_{s,t}$ as a perturbation of $\LL_0$, more precisely we write
\begin{equation}\label{perturb}
\widetilde \LL_{s,t}(\varphi)= \LL_0\biggl (\frac{g_{s,t}}{g_{0,0}} \varphi \biggr )\, .
\end{equation}
Since  $g_{0,0}^{-1}=|f'_0|\in BV_p$ the bound (\ref{algebra})
implies that $\PP_{s,t}(\varphi)= \frac{g_{s,t}}{g_{0,0}} \varphi$
is a bounded operator on $BV_p$. Clearly, $\PP_{0,0}=\id$.
Since $\|g_{s,t}-g_{0,0}\|_{BV_1}\to 0$ 
as $(s,t) \to (0,0)$, the operators $\widetilde \LL_{s,t}$ on $BV_p$ 
have a real positive simple maximal eigenvalue with a spectral
gap, uniformly in  $(s,t)$ close enough to $(0,0)$, by classical perturbation
theory \cite{Ka}.
In particular, the operator $\widetilde \LL_t=\widetilde \LL_{0,t}$ on $BV_p$ has a simple eigenvalue
at $1$, for the fixed point $\tilde \rho_t = \rho_t \circ h_t$,
where $\mu_t=\rho_t \, dx$ is the SRB measure
of $f_t$, and
the rest of its spectrum lies in a disc of strictly smaller radius.
Note that the fixed point of $\widetilde \LL_t^*$ is the measure $\nu_t$ defined by
\begin{equation}\label{stst}
\int \varphi \, dx= \int \varphi \circ h_t\, d\nu_t \, .
\end{equation}
(By definition $\nu_t$ is a probability measure and $\int \tilde \rho_t \, d\nu_t=1$.)

\smallskip

Consider first the case $j=1$.
Lemma~\ref{jensenlemma} below implies  that 
 the map $s \mapsto \PP_{s,t}$ is $C^1$ from 
$\real$ to the Banach space
of $C^1$ maps from $\{  |t|< \epsilon\}$ to bounded operators on $BV_p$,
and
\begin{equation}\label{ppt}
\partial_s \PP_{s,t}|_{s=u}= (\psi \circ h_t)  \PP_{u,t} \,  , \, \, \forall u \in \real .
\end{equation}
Therefore, $s \mapsto \LL_{s,t}$ is $C^1$ from $\real$
to the Banach space of $C^1$ maps from $\{|t|<\epsilon\}$
to bounded operators on $BV_p$, and
$$
\partial_s \widetilde \LL_{s,t} |_{s=u} (\varphi) =\widetilde  \LL_{u,t}( (\psi \circ h_t) \varphi)
\, ,\, \, \forall u \in \real \, .
$$
We are thus in a position to
apply classical perturbation theory of an isolated simple eigenvalue
(see \cite[Ch. VII.1.3]{Ka} for the analytic case, see e.g. \cite[Lemma 3.2]{BH}
for  the differentiable setting). It follows on the one
hand  that,  in a neighbourhood of $(0,0)$, the
maximal
\footnote{Of course, $\log \lambda_{s,t}$ is  the topological pressure
of $\log g_{s,t}$.}
 eigenvalue $\lambda_{s,t}>0$ of $\widetilde \LL_{s,t}$ acting on $BV_p$
is a $C^1$ function of $s$ to the space
of $C^1$ maps from $\{|t|<\epsilon\}$ to $\real$.
On the other hand,
by ``tedious but straightforward calculations" and \cite[Ch. VII.1.5, Ch. II.2.2]{Ka},
(to quote \cite[(5.2)]{KN}),
we have
$$
\partial_{s} (\log \lambda_{s,t})|_{s=0} =
\int \psi \circ h_t \tilde \rho_t \, d\nu_t=  \int \psi \, d \mu_t
$$
(use that 
$\tilde \rho_t$ and $\nu_t$ are the fixed eigenvectors of
$\widetilde \LL_{0,t}$ and its dual).
Since $t\mapsto \partial_{s} (\log \lambda_{s,t})|_{s=0}$ is a $C^1$ function in a
neighbourhood  of zero, we have proved 
Theorem~\ref{pressuretrick} in  the case $j=1$.
If $j\ge 2$, 
apply Lemma~\ref{jensenlemma2} instead of Lemma~\ref{jensenlemma}.
\end{proof}

The following result is the key ingredient in the proof
of Theorem~\ref{pressuretrick}, its proof hinges on Proposition~\ref{glu}
and \cite[Prop. 2.4]{BS1}:

\begin{lemma}\label{jensenlemma}
Let $f_t$ be a $C^{2}$ family of piecewise expanding $C^{3}$ 
unimodal maps in the topological class of $f_0$.
For any $p> 1$  there exists $\epsilon_p>0$ so that
for any $\psi: I \to \real$ which is $C^{1+\Lip}$,
the map $ s\mapsto g_{s,t }$ defined by (\ref{gst})
is $C^1$ from  $\real$ to the Banach space
of $C^1$ maps from $\{ |t|<\epsilon_p\}$ to $BV_p$.
In addition, recalling the notation (\ref{starr}),
\begin{equation}\label{easy}
\partial_s g_{s,t}|_{s=u}= (\psi \circ h_t) \, g_{u,t} \, , \, \, \forall u \in \real \, .
\end{equation}
\end{lemma}

In fact, $s$-analyticity holds in Lemma~\ref{jensenlemma}, but we shall not need this.

\begin{proof} [Proof of Lemma ~ \ref{jensenlemma}]
Fix $p >1$. For every $x \ne c$, all small $t$, and all $s_1 <s_2$
in $\real$, there exists $s_3\in [s_1, s_2]$ so that
\begin{equation}\label{starttaylor}
g_{s_1,t}(x)- g_{ s_2, t}(x) - (s_1- s_2) \psi(h_t(x)) g_{ s_2, t}(x)
=(s_1- s_2)^2 (\psi(h_t(x))) ^2 g_{ s_3, t}(x) \, .
\end{equation}
(Just use the Taylor formula for $s \mapsto g_{s,t}(x)$ and the intermediate value
theorem.)

So, to prove both differentiability and (\ref{easy}), it suffices to see that the three maps
$$
t \mapsto g_{s,t} \, ,\, t \mapsto (\psi\circ h_t) g_{s,t}\, , \, 
t \mapsto(\psi\circ h_t)^2 g_{s,t} \, ,
$$
are $C^1$ from a neighbourhood of $0$ to $BV_p$, uniformly in $s$
in any compact set $K\subset \real$.

In view of this, we first study the
maps $t \mapsto h_t(x)$. 
By \cite[Proposition 2.4]{BS1}, there
exists $\tilde \epsilon >0$ so that the set of maps
$\{    t\mapsto h_t(x), \ x \in I    \}$
is  bounded in $C^{1+\Lip}([-\tilde \epsilon,\tilde \epsilon])$.  Differentiating  with respect to $t$ the  
equation $h_t \circ f_0 = f_t \circ h_t$, and setting
$\alpha_t = \partial_t h_t \circ h_t^{-1}$, we get
$$
\alpha_t(f_t(c))=\partial_t f_t (c)\, ,
\quad \partial_t f_t (x) = \alpha_t(f_t(x))-f'_t(x) \alpha_t (x)  \, ,\forall x \ne c
\, , |x|< \tilde \epsilon \, .
$$
Since $\alpha_t(c)=0$ this implies
$J(f_t,\partial_t f_t)=0$ for  $|t|<\tilde \epsilon$
(recall (\ref{bdd})), 
so, for any fixed
$$ 
\beta \in (1/p, 1) \, 
\mbox{ (we may and shall assume also that $\beta <1/\sqrt{p}$),}
$$
Proposition~ \ref{glu} gives  $C$ and $\epsilon_p>0$ so that 
\begin{equation}\label{fromglu}
|\alpha_t|_\beta \le C  \, , \quad \forall |t|<\epsilon_p \, .
\end{equation}
Let $\alpha^{\eta}_t$ be the $\eta$-regularisation (in the variable $x$)
of $\alpha_t$, that is the convolution $\alpha^\eta_t(x)=
\int \alpha_t (y) \kappa_\eta(x-y)\, dy$ of $\alpha_t$
with a convolution kernel  $\kappa_\eta(x)=\eta^{-1}\kappa(x/\eta)$, where
the $C^\infty$ function $\kappa: \real\to \real_+$ is supported in $[-1,1]$,
and $\int \kappa(x)\, dx=1$.
Note for further use that 
(\ref{fromglu}) gives $\widetilde C$ so that, for all  $|t|<\epsilon_p$,
\begin{equation}\label{convol}
\| \alpha^{\eta}_t\|_{C^{1+1/p}} \le  \frac{\widetilde C}{\eta^{1+1/p-\beta}}\, ,\, 
\| \alpha^{\eta}_t\|_{\beta} \le  \widetilde C
\, , \,  \, | \alpha^{\eta}_t-\alpha_t |_{1/p}  \le \widetilde C 
\eta^{\beta-1/p}\, , \, \,   \forall \eta\in (0,1)\, .
\end{equation}

We now consider $t \mapsto g_{s,t}$.
For $x\ne c$, we have
\begin{equation}\label{first}
\partial_t g_{s,t}(x)=e^{s \psi(h_t(x))}
\left [ \frac
{ \psi'(h_t(x)) \alpha_t(h_t(x))}{|f'_t(h_t(x))|}
-\frac{\partial_t\bigl (|f'_t(h_t(x))| \bigr )  }
{|f'_t(h_t(x))|^2}  \right ]\, ,
\end{equation}
where
\begin{equation}\label{second}
\partial_t \bigl (|f'_t(h_t(x))|\bigr )=-\sgn(x) \bigl ( f''_t(h_t(x)) \alpha_t(h_t(x))+
\partial_t f'_t(h_t(x)) \bigr )\, .
\end{equation}
We claim that the 
function $x \mapsto \partial_t g_{s,t}(x)$ 
has bounded $BV_{1/\beta}$ norm,  uniformly in $s\in K$ and $|t|<\epsilon_p$.
Indeed, decomposing
$$\partial_t g_{s,t}= b_{s,t} \circ h_t\, , $$ 
note that  each $h_t:I\to I$ 
is a homeomorphism leaving both $[-1,c]$ and $[c,1]$ invariant,
while $b_{s,t}$ is $\beta$-H\"older
on $[-1,c)$ and $(c,1]$, uniformly in $s\in K$
and $|t|<\epsilon_p$ (because $\psi'$ is $C^\beta$,
$f_t$ is a $C^2$ family of $C^3$ maps\footnote{This implies in particular
that $x \mapsto \partial_t f_t$ is $C^2$ in $x$, uniformly in $t$
and $\partial_x \partial_t f_t=\partial_t f'_t$.}, and
$\alpha_t$ is $\beta$-H\"older, uniformly
in  $|t|<\epsilon_p$),
and $\sup_{s\in K,|t|<\epsilon_p}|b_{s,t}(c_+)- b_{s,t}(c_-) |< \infty$
(using $\sup_{|t|<\epsilon_p} \|f_t\|_{\BB^{2+\beta}}<\infty$).

To conclude, it suffices to prove that our candidate
$ b_{s,t} \circ h_t\in BV_p$ is really the $t$-derivative of $g_{s,t}$
(uniformly in $s$), that is,
\begin{equation}\label{aa}
\lim_{t_2 \to t_1} \sup_{s \in K}
\biggl \| \frac{g_{s,t_2} - g_{s,t_1}}{t_2-t_1}-  b_{s,t_1} \circ h_{t_1}\biggr \|_{BV_p} = 0 \, ,
\forall |t_1|<\epsilon_p\, ,
\end{equation}
and that this derivative is continuous in $t$ (uniformly in $s$), that is,
\begin{equation}\label{bb}
\lim_{\tilde t \to t_1} \sup_{t_2 \in [t_1,\tilde t]} \sup_{s \in K}
\| b_{s,t_2} \circ h_{t_2}- b_{s,t_1} \circ h_{t_1}\|_{BV_p} = 0 \, ,
\forall |t_1|<\epsilon_p \, .
\end{equation}
We first prove (\ref{bb}). 
Decomposing 
\begin{equation}\label{bothterms}
b_{s,t_2} \circ h_{t_2}- b_{s,t_1} \circ h_{t_1}=
(b_{s,t_2} - b_{s,t_1} )\circ h_{t_2}
+b_{s,t_1} \circ h_{t_2}- b_{s,t_1} \circ h_{t_1}\, ,
\end{equation}
we focus first on the second term in the  right-hand-side of (\ref{bothterms}).
Let $\delta >0$ be such that $f'_t$, $f''_t$ and $\partial_t  f'_t$
restricted to $[-1,c]$ and $[c,1]$, respectively,
extend to $C^1$ functions of $x$ on $[-1-\delta,c+\delta]$
and $[c-\delta, 1+\delta]$, respectively, for all $|t|< \epsilon_p$.
Denote by $b^{\eta,-}_{s,t}$ the function obtained from
$b_{s,t}$ by substituting $\alpha_t$ with $\alpha^\eta_t$, and also
$\psi'$,  and the extensions to $[-1-\delta,c+\delta]$
of $f''_t|_{[-1,c]}$,  $\partial_t  f'_t|_{[-1,c]}$, 
with their $x$-convolutions with $\kappa_\eta$,
for small $\eta >0$ (to be determined later).
Define $b^{\eta,+}_{s,t}$ similarly, using $[c-\delta, 1+\delta]$, and
set $b^{\eta}_{s,t}(x)=b^{\eta,+}_{s,t}(x)$
if $x>c$ and $=b^{\eta,-}_{s,t}(x)$ if $x <c$. 
Since $\beta <1$ and $\psi'$ is Lipschitz, it is easy to see that
there exists $\widehat C>0$ so that
for all $\eta \in (0,1)$  
\begin{equation*}\label{convol2}
\max \bigl (\sup_{s \in K}|(b^\eta_{s,t_1}|_{(-\infty,c)})'|_{1/p} ,
\sup_{s \in K}|(b^\eta_{s,t_1}|_{(c,\infty)})'|_{1/p} \bigr )
\le \frac{\widehat C}{ \eta^{1+1/p-\beta}}\, , \, \,  \forall |t_1|<\epsilon_p\, 
\, .
\end{equation*}
(Use  the first two estimates of (\ref{convol}), and the analogous
bounds for the regularisations
of $\psi'$ and  $f''_t$, $\partial_t f'_t$.)
Therefore, by the fundamental theorem of
calculus and the H\"older (or Jensen) inequality, there exists $\bar C>0$
so that for all $s\in K$, all $|t_1|<\epsilon_p$, $|t_2|<\epsilon_p$, all  $\eta\in (0,1)$,
and any $x_0 < x_1 < \cdots < x_N\le c$,
\begin{align}
\nonumber &\sum_{i=0}^{N-1} |b^\eta_{s,t_1} ( h_{t_2}(x_i))- b^\eta_{s,t_1} ( h_{t_1}(x_i))
-b^\eta_{s,t_1} ( h_{t_2}(x_{i+1}))+b^\eta_{s,t_1} ( h_{t_1}(x_{i+1}))|^p \\
\nonumber &=
\sum_i \bigl |
\int_{t_1}^{t_2} \partial_t (b^\eta_{s,t_1} ( h_{t}(x_i)))\, dt
- \int_{t_1}^{t_2}\partial_t (b^\eta _{s,t_1} ( h_{t}(x_{i+1})))\, dt \bigr |^p \\
\nonumber &\le \sum_i \int_{t_1}^{t_2}
|  (b^\eta_{s,t_1} )'( h_{t}(x_i))\alpha_t(h_t(x_i))- (b^\eta _{s,t_1}) '( h_{t}(x_{i+1}))
\alpha_t(h_t(x_{i+1}))|^p\, dt\\
\nonumber&=\int_{t_1}^{t_2} \sum_i 
| (b^\eta_{s,t_1})' ( h_{t}(x_i))\alpha(h_t(x_i))-
 ( b^\eta _{s,t_1})' ( h_{t}(x_{i+1}))\alpha_t(h_t(x_{i+1}))|^p\, dt
\\
 \nonumber &\le  |t_2-t_1| 
\bigl ( \sup_{x \in (-\infty,c)} |(b^\eta_{s,t_1})'(x) | \sup_t \| \alpha_t\circ h_t \|_{BV_p}
+\sup_{x,t} |\alpha_t| |(b^\eta_{s,t_1}|_{(-\infty,c)})'|_{1/p} \bigr ) \\
\label{okk}
&\le   \bar C \frac{ |t_2-t_1| }{\eta^{1+1/p-\beta}}   \, .
\end{align}
(We used (\ref{fromglu}) in the last inequality.)
The same bounds hold for  $c \le x_0 < x_1 < \cdots < x_N$, and it is easy
to estimate the jump of
$b^\eta_{s,t_1}\circ h_{t_2} -b^\eta_{s,t_1}\circ h_{t_1}$
at $x=c$ uniformly in $s$ and $t_1$, $t_2$.

We next analyse the contribution of $b_{s,t_1}-b^\eta_{s,t_1}$
to the second term of (\ref{bothterms}). For this,
observe   that if $h$ is an orientation preserving homeomorphism fixing
$c$ and $b$ is $\beta$-H\"older on
$[-\infty,c]$ and $[c,\infty]$, then
$\|b \circ h\|_{BV_{p}} \le |b|_{(-\infty,c)}|_{\beta}
+|b|_{[c,\infty)}|_{\beta}+|b(c_+)-b(c_-)|$.
Then,   
the last bound of (\ref{convol}) and its analogue for the
$\eta$-regularisation of $\psi'$, $f''_t$ and $\partial_t f'_t$ give a constant
$C'$ so that for all $|t_1|<\epsilon_p$, $|t_2|<\epsilon_p$
and $\eta\in (0,1)$
\begin{align}\nonumber
\sup_{s\in K} \|(b_{s,t_1}-b^\eta_{s,t_1}) \circ h_{t_2}-
(b_{s,t_1}- b^\eta_{s,t_1} )\circ h_{t_1}\|_{BV_p}
&\le 2 \sup_{s,t} 
\|(b_{s,t_1}-b^\eta_{s,t_1}) \circ h_{t}\|_{BV_p}\\
\label{akk}&\le  C' \eta^{\beta-1/p}  \, .
\end{align}
Taking  $\xi \in (0,1)$ and setting
$
\eta= (t_2-t_1)^{\frac{\xi}{1+1/p-\beta}} 
$,
we get from (\ref{okk}--\ref{akk}) that
$\lim_{\tilde t\to t_1}\sup_{t_2\in [t_1, \tilde t]}\sup_{s\in K}\|b_{s,t_1} \circ h_{t_2}- b_{s,t_1} \circ h_{t_1}\|_{BV_p}=0$.

To analyse the first term of (\ref{bothterms}), we start by noticing that since
$t \mapsto \partial_t h_t$ is Lipschitz, there exists a 
set $\DD_p\subset (-\epsilon_p,\epsilon_p)$ of full Lebesgue measure
so that $\partial_t h_t$ is 
differentiable at all $t$ in  $\DD_p$. 
Differentiating   twice  $ f_t \circ h_t(x)=h_t \circ f(x)$
with respect to $t$  and
\footnote{This is  similar the proof of \cite[Proposition 2.4]{BS1}, but we will
make a more careful analysis of what was called $F_i$ there.}
setting $\alpha^2_t= \partial^2_{tt} h_t\circ h_t^{-1}$,
we obtain for all $x \ne c$ and all $t \in \DD_p$ that
\begin{equation}\label{j=2}
 f''_t   (x) \alpha_t(x)^2 +  2\partial_{t} f'_t (x) \alpha_t(x) + \partial_{tt} f_t(x)
= \alpha^2_t(f_t(x))- f'_t(x)  \alpha^2_t(x) 
  \, .
\end{equation}
The left-hand-side of the above TCE
is $\beta$-H\"older in $[-1,c]$ and $[c,1]$ and continuous in $I$, 
since $\alpha_t(c)=0$ for every small $t$, so it is $\beta$-H\"older continuous. 
Therefore, by 
Proposition~  \ref{glu}, there exist $\epsilon_p>0$ and  a constant $C''$ so that 
\begin{equation}\label{DDp}
|\alpha^2_t|_\beta \le C''\, ,\quad \forall t \in \DD_p\, .
\end{equation}
The fundamental theorem of calculus holds for the Lipschitz
(and therefore almost everywhere differentiable) function
$t\mapsto b_{s,t}$ and gives
\begin{equation}\label{fund}
(b_{s,t_2}-b_{s,t_1}) h_{t_2}(x)=
\int_{t_1}^{t_2} \partial_ t b_{s,t}( h_{t_2}(x))\, dt
\, , \forall x \ne c\, .
\end{equation}
The first term of (\ref{bothterms})
may then be estimated via the H\"older inequality
and  the fundamental theorem of calculus
(\ref{fund}), as in (\ref{okk}), but exploiting (\ref{DDp}) 
instead of  using $\eta$-regularisation. Details are left to the reader.

Finally, to show (\ref{aa}), start from
$$
g_{s,\tilde t}(x) - g_{s,t_1}(x)-  (\tilde t-t_1) b_{s,t_1} (h_{t_1}(x))
=\int_{t_1}^{\tilde t}( b_{s,t_2} ( h_{t_2} (x)) -  b_{s,t_1} (h_{t_1}(x)))\, dt_2\, ,
$$ 
for all $x\ne c$, and use the H\"older inequality
and (\ref{bb}) (details are left to the reader).

The analysis of the maps $t \mapsto (\psi\circ h_t) g_{s,t}$ and
$t \mapsto(\psi\circ h_t)^2 g_{s,t}$ goes along exactly the same lines.
\end{proof}

For the higher regularity statement in Theorem ~\ref{pressuretrick}, we use
the following result
(again, analyticity in $s$ holds):

\begin{lemma}\label{jensenlemma2}
Let $j\ge 2$. Let $f_t$ be a $C^{j+1}$ family of piecewise expanding $C^{j+2}$ 
unimodal maps in the topological class of $f_0$.
For any $p> 1$  there exists $\epsilon_p>0$ so that
for any $\psi: I \to \real$ which is $C^{j+\Lip}$,
the map $ s\mapsto g_{s,t }$ defined by (\ref{gst}) is $C^1$ from  $\real$ to the space 
of $C^j$ maps from  $\{ |t|<\epsilon_p\}$ to $BV_p$,
and, recalling (\ref{starr}),
$\partial_s g_{s,t}|_{s=u}= (\psi \circ h_t) \, g_{u,t}$.
\end{lemma}

\begin{proof}
Since the family $f_t$ is $C^{j+1}$, the set
$\{    t\mapsto h_t(x), \ x \in I    \}$
is  bounded in $C^{j+\Lip}$ by \cite[Proposition 2.4]{BS1}.  
Let $\beta \in (1/p, 1)$ (with $\beta <1/\sqrt p$, say). Assume first $j=2$. 
Then, by (\ref{j=2}), the function $\alpha^2_t= \partial^2_{tt} h_t\circ h_t^{-1}$,
is well-defined for all $|t|< \epsilon_p$ and there exists $C$ so that 
$|\alpha^2_t|_\beta \le C$ for every  $|t|<\epsilon_p$.
For $j\ge 3$, a higher order TCE similar to
(\ref{j=2}) gives that $\alpha^j_t= \partial^j_{t^j} h_t\circ h_t^{-1}$
is $\beta$-H\"older for all $|t|< \epsilon_p$.
 We put $\alpha^1_t=\alpha_t$.

Then, computing $\partial^j_{t^j} g_{s,t}(x)$
at $x\ne c$ gives $ b^{(j)}_{s,t} ( h_t(x))$, where $b^{(j)}_{s,t}$
is an expression  involving derivatives of order at most $j$
of $\psi(x)$,  functions
$\alpha_t^\ell$,  for $1\le \ell \le j$, and derivatives
(in $x$,  $t$, or mixed)  of total
order at most $j$ of $f'_t(x)$,  in the numerator,
and $|f'_t(x)|^m$ for $m\ge 1$ in the denominator.
Our differentiability assumptions on $\psi$ and the family $f_t$ then allow us to proceed
as in the proof of Lemma~\ref{jensenlemma} (using Taylor series
of higher order).
\end{proof}

\section{Recovering the linear response formula}
\label{secc}

Here we give a slightly different proof of the differentiability of $R(t)=\int \psi\, d\mu_t$, 
where $\mu_t$ is the SRB measure of $f_t$,
still relying heavily on Proposition~ \ref{glu} (via
Lemma~\ref{gla}).
The advantage with respect to Theorem~\ref{pressuretrick}  is that
we recover the formula for $\partial_t R(t)|_{t=0}$, and we need only assume that
$\psi$ is $C^0$. (In particular, this gives a new proof of \cite[Theorem 5.1]{BS1}.)
We also get new information in Corollary~\ref{thebest} by combining
Theorems~\ref{pressuretrick} and \ref{theformula}.

We need  notation.
By \cite[Proposition 3.3]{Ba},
we may decompose the invariant density
of a piecewise expanding $C^3$ unimodal mixing map $f_t$
as $\rho_t=\rho_{reg,t}+\rho_{sal,t}$,
where $\rho_{reg,t}\in BV \cap C^0$,  $\rho'_{reg,t}\in BV$, and
$$
\rho_{sal,t}= \sum_{k=1}^{M_f}  s_{k,t} H_{c_{k,t}} \, .
$$
(Here, $H_u(x)$ denotes the Heaviside function $H_u(x)=-1$ if
$x<u$, $H_u(x)=0$ if $x>u$ and $H_u(u)=-1/2$.)
If $M_f=\infty$ then it is not difficult to show that (see e.g. 
\cite{Ba, BS1}, noting that if $c_{1,t}$ is preperiodic
but not periodic our notation is slightly different than the
notation there) 
\begin{equation}\label{decay}
s_{k,t}=\frac{s_{1,t}}{(f^{k-1})'(c_{1,t})} \, ,\quad \forall k\ge 1\, .
\end{equation}
We simply write  $\rho_0=\rho=\rho_{reg}+\rho_{sal}$.

To compute the formula for the derivative,
we shall assume, as in \cite{BS1},  that  $v= \partial_t f_t |_{t=0}$ is of the
form $v=X \circ f_0$ for a $C^2$ function $X:I \to \real$. 
\footnote{See also the beginning of \cite[Section 17]{RuStr}.}

\begin{theorem}\label{theformula}
Let $f_t$ be a $C^2$ family of piecewise expanding $C^3$ unimodal maps.
Assume that $f_0$ is good and mixing, that
$f_t$ is tangent to the topological
class of  $f_0$, and that
$v=\partial_t f_t |_{t=0}=X \circ f_0$ for a $C^2$ function $X$. Then, as Radon measures,
\begin{equation}\label{formula}
\lim_{t \to 0} \frac{\mu_t-\mu_0}{t}=- \alpha  \rho'_{sal}- (\id- \LL_0)^{-1} (X' \rho_{sal}+(X\rho_{reg})') \, dx \, ,
\end{equation}
where the function $\alpha$ is given by (\ref{bdd}), and the
operator $\LL_0=\widetilde \LL_{0,0}$ is defined
by (\ref{xfer1}).
In addition, $\alpha$  is $\beta$-H\"older
for any $\beta <1$.
\end{theorem}

\begin{proof}
Set $f=f_0$ for convenience.
By Proposition~\ref{KLbd}, we can assume that $f_t$ lies
in the topological class of $f$, denoting the conjugacies by $h_t$ as usual.
The transfer operator $\widetilde \LL_t=\widetilde  \LL_{0,t}$ for $f$ and the
weight $|f'_t \circ h_t|^{-1}$ (recall (\ref{xfer1})) is conjugated to the transfer operator for $f_t$ and $|f'_t|^{-1}$
defined by
$$\LL_t \varphi(x)=\sum_{f_t(y)=x}\frac{\varphi(y)}{|f'_t(y)|}\, ,
$$ 
via
\begin{equation}
\label{zero}
\widetilde \LL_t (\varphi \circ h_t)= \LL_t(\varphi) \circ h_t \, .
\end{equation}
Since $h_t$ is a homeomorphism, it gives rise to an isometry of $BV_p$, and
(\ref{zero}) together with the  main Theorem of \cite{W} applied  to $\LL_t$
imply that there exist $p_0>1$ so that for any 
$p\in [1, p_0)$ there exists $\epsilon_p>0$ so that for all $|t|< \epsilon_p$
the operator
$\widetilde \LL_t$ acting on $BV_p$ has a maximal eigenvalue
equal to $1$, which is simple, and  the rest of the
spectrum lies in a disc of strictly smaller radius
(i.e.,  $\widetilde \LL_t$ has a spectral gap). 
The fixed points
of $\widetilde \LL_t$ and its dual, $\widetilde \rho_t=\rho_t \circ h_t$ and
$\nu_t$ from (\ref{stst}), were  introduced in the proof of Theorem ~ \ref{pressuretrick}.
We can alternatively see that $\tilde \LL_t$ has a spectral gap on $BV_p$
by noting that it is a small multiplicative perturbation of $\LL_0=\tilde \LL_0$ on $BV_p$:
Recall that $\lim_{t \to 0}\|g_{0,t} -g_{0,0}\|_{BV_1}=0$ and(\ref{perturb}). Observe for further use that this implies 
$\lim_{t \to 0} \|\tilde \rho_t - \rho_0\|_{BV}=0$.

From now on, we fix  $p\in (1,p_0)$.

We next show that  $t\mapsto \tilde \rho(t)\in BV_p$ and $t \mapsto \nu_t\in BV_p^*$ 
are differentiable at $t=0$. 
By \cite[Prop. 2.4, Cor. 2.6]{BS1} $v$ is horizontal
for $f_0$, $t \mapsto h_t(x)$ is differentiable, uniformly
in $x \in I$, and $\alpha=\partial_t h_t|_{t=0}$
is continuous, with $\alpha(c)=0$, $\alpha(c_1)=X(c)$, and
$\alpha$ is the unique bounded solution (\ref{bdd})
to the TCE (\ref{tceeq}).  In addition,
Proposition ~\ref{glu} gives that $\alpha$
is $\beta$-H\"older for arbitrary $\beta <1$
(we shall take $\beta \in (1/p, 1/\sqrt p)$).

Our assumptions on $f_t$ then imply that $v'$ is $C^1$ and
the following operator is bounded on $BV_p$:
\begin{equation}\label{deriv0}
\MM \varphi (x)=- \sum_{f(y)=x} \frac{f''(y) \alpha(y) + v'(y)}{|f'(y)|f'(y)} \varphi(y) \, .
\end{equation}
(Write $\MM$ as $\LL_0$ composed with a multiplication
operator, like in (\ref{perturb}) and use (\ref{algebra}).)
Lemma~\ref{gla} below easily implies that $t \mapsto \tilde \LL_t$ is differentiable
as an operator on $BV_p$, and that
\begin{equation}\label{deriv}
\partial_t \tilde \LL_t |_{t=0}=\MM \, .
\end{equation} 
As in the proof of Theorem~ \ref{pressuretrick}, 
perturbation theory then  gives that
$t \mapsto \tilde \rho_t \in BV_p$ and $t \mapsto \nu_t\in BV_p^*$ 
are differentiable at $t=0$.

\smallskip

We next show that $t \mapsto \mu_t= \rho_t\, dx$ is differentiable as a Radon measure,
exploiting  the formula for $\MM$ to get the claimed formula
for $\partial_t \mu_t|_{t=0}$. Fix $\psi : I \to \complex$  continuous.
Since $\tilde \rho_0=\rho_0$, we can decompose
\begin{align}\label{deccomp}
\int \psi \rho_t \, dx - \int \psi  \rho_0 \, dx&=
\int \psi \rho_t \, dx - \int \psi \tilde \rho_t \, dx +\int \psi \tilde \rho_t \, dx - \int \psi  \tilde \rho_0 \, dx\, .
\end{align}
We shall now see that 
\begin{equation}\label{otherterm}
\lim_{t\to 0}\frac {\int \psi (\rho_t - \rho_t \circ h_t) \, dx}{t}
= -\int \psi \alpha \rho_0' \, .
\end{equation}
In view of (\ref{otherterm}), note first that $s_{k,t}\to s_k$ as $t\to 0$: Indeed, 
since $\lim_{t \to 0} \|\tilde \rho_t - \rho_0\|_{BV}=0$
we have
(in $BV$)
\begin{equation}\label{important}
\lim_{t \to 0}\tilde \rho_t =
\lim_{t \to 0}\biggl (\rho_{reg,t} \circ h_t + \sum_{k=1}^{M_f}  s_{k,t} H_{c_{k}}
\biggr )= \tilde \rho_0=
\rho_{reg}+\sum_{k=1}^{M_f}  s_{k} H_{c_{k}}  \, .
\end{equation}
(We gave another proof of 
$\lim_{t \to 0}s_{k,t}= s_k$ in Step 1 of \cite[Proof of Theorem 5.1]{BS1}.)

Decompose $\rho_t-\rho_t \circ h_t$ in (\ref{otherterm})
into  $\rho_{sal,t}-\rho_{sal,t} \circ h_t + \rho_{reg,t}-\rho_{reg,t} \circ h_t$.
For the singular term, we have
 in the sense of Radon measures:
\begin{equation}\label{a}
\lim_{t \to 0}\frac{\rho_{sal,t }- \rho_{sal,t} \circ h_t}{t}=
 -\sum_{k =1} ^{M_f}  \alpha(c_k) s_k \mbox{Dirac}_{c_k} =
-\alpha \rho'_{sal} \, . 
\end{equation}
(Just use that $s_{k,t}\to s_k$ and (\ref{decay}), which implies
the $s_{k,t}$ decay exponentially
in $k$ uniformly in $t$.\footnote{Claim (\ref{a})
 was also proved in Step 1 of \cite[Proof of Theorem 5.1]{BS1}.})

We claim that the contribution of the
regular term $\rho_{reg,t}$  in the decomposition of (\ref{otherterm})
is $-\int \psi \alpha \rho_{reg}'\, dx$. Indeed, 
recall that $\rho_{reg,t}' \in BV$. In  fact, we may decompose
\begin{equation}\label{dec0}
\rho_{reg,t}=\rho_{regreg,t}+ \rho_{regsal,t}\, ,
\end{equation}
with (see Step~ 2 of  \cite[Proof of Theorem 5.1]{BS1}) $\rho_{regreg,t}$ is $C^1$ 
(uniformly in $t$) and
$$
 \rho'_{regsal,t}
=\sum_{k=1}^{M_f}  s'_{k,t} H_{c_{k,t}}\, ,
$$
where the $s'_{k,t}$ decay exponentially uniformly in $t$.
In fact (see \cite[(69), (71)]{BS1})
\footnote{If $c$ is periodic then $(\rho_{reg,t})'(c)$ may be undefined,
but $(\rho_{reg,t})'(c_\pm)$ are both defined.} 
\begin{equation*}
s'_{k,t}=E'_{k,t}-E_{k,t}\, , \mbox{ with }
E'_{k,t}=\frac{s'_{k-1,t}}{(f'_t(c_{k-1,t})^2)} \,, \, k\ge 2\, , 
\end{equation*}
and 
\begin{align*}
E_{k,t} &=
\frac{s_{k-1,t} 
f''_t(c_{k-1,t})}{(f'_t(c_{k-1,t}))^3} \, , \,\, k\ge 2\, , \quad
 E'_{1,t}=-\frac{(\rho_{reg,t})' (c)}{(f'_t(c_-))^2}
+ \frac{(\rho_{reg,t})' (c)}{(f_t'(c_+))^2}\\
\nonumber  E_{1,t}&= \biggl (-
 \frac{\rho_{reg,t} (c)f''_t(c_-)}{(f'_t(c_-))^3}
+  \frac{\rho_{reg,t} (c)f''_t(c_+)}{(f'_t(c_+))^3}\biggr ) \\
\nonumber
&\qquad\qquad\qquad\qquad+\sum_{k \ge 2, c_{k-1,t}>c} s_{k-1,t}
\biggl (\frac{f''_t (c_-)}{(f'_t(c_-))^3} -\frac{f''_t (c_+)}{(f'_t(c_+))^3} \biggr )\, .
\end{align*}

The argument above giving $s_{k,t}\to s_k$ also yields $s'_{k,t}\to s'_k$
(just differentiate once). 
Therefore, just like in (\ref{a}), we have
\begin{align}\label{b}
\lim_{t \to 0}\int \psi \frac{\rho_{regsal,t }- \rho_{regsal,t} \circ h_t}{t} \, dx
&= -
\sum_{k =1} ^{M_f} \int  \alpha(x) s'_k \psi (x) H_{c_k}(x)\, dx \\
\nonumber &=\int \psi \alpha \rho_{regsal}' dx\, .
\end{align}

In view of handling the term $\rho_{regreg,t}$ from
(\ref{dec0}), observe that 
\begin{equation}\label{this}
\lim_{t\to 0} \|\varphi - \varphi \circ h_t\|_{BV}= 0 \, ,
\quad \forall \varphi \in C^1 \, .
\end{equation}
Indeed, for $\delta >0$ and any partition $x_0 < \ldots x_i < x_{i+1}< \ldots < x_n$ 
let 
$N\le n$ be so that
$\min(x_N,\inf_t h_t(x_N)) > 1-\delta$, 
and since $|h_t(y)-y|=O(t)$ uniformly in $y$, take $t_0$
so that $|x_i-h_t(x_i)| < \delta/N$ for all $i\le N$ and $|t|<t_0$.
Then  use 
\begin{align*}
&\sum_{i=0}^{n-1} |\varphi(x_i)-\varphi(x_{i+1})
-\varphi(h_t(x_i))+\varphi(h_t(x_{i+1}))|\le2\sum_{i=0}^N |\varphi(x_i)-\varphi(h_t(x_{i}))| \\
&\, +\sum_{i=N+1}^{n-1} |\varphi(x_i)-\varphi(x_{i+1})| + 
\sum_{i=N+1}^{n-1} |\varphi(h_t(x_i))-\varphi(h_t(x_{i+1}))|
\, .
\end{align*}
Since the $C^1$ norm of $\rho_{regreg,t}$
is bounded uniformly in $t$,
(\ref{this}) easily implies that 
\begin{equation}\label{A}
\lim_{t \to 0} \|\rho_{regreg,t}- \rho_{regreg}\|_{BV}= 0 \, .
\end{equation}
(Note for the record that,
since  $\|\tilde \rho_t - \rho_0\|_{BV}=0$ and $s_{k,t}\to s_k$,
 $s'_{k,t}\to s'_k$, with $t$-uniformly exponentially
decaying $s_{k,t}$, $s'_{k,t}$,   this implies $\lim_{t \to 0} \|\rho_{reg,t}- \rho_{reg}\|_{BV}= 0$.
\footnote{$\lim_{t \to 0} \|\rho_{reg,t}- \rho_{reg}\|_{BV}= 0$
can also be proved from Keller-Liverani techniques on the spaces
in \cite{BS1}, without using the fact that $\tilde \rho_t \to \rho_0$ in the
$BV$ norm.  This implies $s'_{1,t}\to s'_1$, like we proved $s_{1,t}\to s_1$ in
\cite[Proof of Theorem 5.1]{BS1}.})
Then, by the mean value theorem
and the $x$-uniform differentiability of $t\mapsto h_t(x)$
\begin{align}
\nonumber &\lim_{t\to 0}\int \psi \frac{\rho_{regreg,t }- \rho_{regreg,t} \circ h_t}{t}  dx\\
\nonumber &\qquad \qquad
= \lim_{t\to 0}\int \psi(x) \frac{\rho_{regreg,t }(x)- \rho_{regreg,t} ( h_t(x))}{x-h_t(x)}
\frac{x-h_t(x)}{t}  dx\\
\nonumber &\qquad\qquad 
= \lim_{t\to 0}\int \psi (x)\rho'_{regreg,t }(x_t) \frac{x-h_t(x)}{t} \, dx\\
\nonumber &\qquad\qquad
=- \lim_{t \to 0}\int \psi (x)\rho'_{regreg,t }(x_t) \alpha(x) \, dx \\
\label{c} &\qquad\qquad =- \int \psi(x) \rho'_{regreg,0} (x)\alpha(x) \, dx \, ,
\end{align}
where $x_t$ is in the interval between $x$ and $h_t(x)$,
and we used in the last line
that  $\rho_{regreg,t}'$ is continuous on the compact interval $I$, uniformly in $t$,
together with (\ref{A}).
Putting (\ref{a}-\ref{b}-\ref{c}) together, we find (\ref{otherterm}).

We now turn to the estimation
of the term $(\int \psi \tilde \rho_t \, dx - \int \psi  \tilde \rho_0 \, dx)/t$ 
from (\ref{deccomp}).
In view of this, note that
(\ref{deriv}) implies that  (as operators on $BV_p$)
$$
\partial_t (z-\tilde \LL_t)^{-1}|_{t=0}=
(z-\LL_0)^{-1} \MM (z-\LL_0)^{-1} \, .
$$
Therefore, writing the spectral projectors as Cauchy integrals, we get by a
simple residue computation, since
$ (z-\LL_0)^{-1} \rho_0=\rho_0/(z-1)$,
\begin{equation}\label{resid}
\partial_t (\nu_t(\rho_0) \tilde \rho_t)|_{t=0} =  (\id-\LL_0)^{-1}(\id-\Pi_0) \MM \rho_0\, ,
\end{equation}
where $\Pi_0(\varphi)=\rho_0 \int \varphi\, dx$.

Next, we claim that we have
(in $BV_p$)
\begin{align}\label{ccclaim}
-\alpha \rho_{reg}' 
+  (\id-\LL_0)^{-1}(\id-\Pi_0)( \MM \rho_0)
 = -(\id-\LL_0)^{-1}(\id-\Pi_0)(X'\rho_0 +X \rho_{reg}') \, .
\end{align}
(Recall from \cite[Proof of Proposition~4.4]{Ba} that $\Pi_0 (X'\rho_0+X \rho_{reg}')=0$.)
Since the TCE (\ref{tceeq}) implies, using $v'=(X'\circ f )\cdot f'$,
\begin{align*}
 \MM \rho_0(x)= (X(x)-\alpha(x))\sum_{f(y)=x} \frac{f''(y) }{|f'(y)|f'(y)^2}\rho_0(y)
 - X'(x) \rho_0(x) \, ,
\end{align*}
to prove (\ref{cclaim}), it suffices to show
\begin{align*}
-\alpha \rho_{reg}' 
+  (\id-\LL_0)^{-1}(\id-\Pi_0)(X- \alpha ) ( \widetilde \MM \rho_0)
 =- (\id-\LL_0)^{-1}(\id-\Pi_0)(X \rho_{reg}')\, ,
\end{align*} 
where 
$\widetilde \MM \varphi(x)= \sum_{f(y)=x} \frac{f''(y) }{|f'(y)|f'(y)^2}\varphi(y)$.
Note that \footnote{This is \cite[(70)]{BS1}.}  if $x\in [-1,c_1)$ 
is not along the postcritical orbit we have,
using $(\rho_{reg})'(y)=(\rho_0)'(y)$ if $y$ is not on the postcritical orbit,
\begin{equation}\label{star}
(\rho_{reg})'(x)=(\rho_0)'(x)=(\LL_0(\rho_0))'(x)
=\sum_{f(y)=x} \frac{(\rho_{reg})'(y)}{|f'(y)|f'(y)}-
\frac{\rho_0(y) f''(y)}{|f'(y)| (f'(y))^2}\, .
\end{equation}
Therefore, for any $x\in I$ which is not on the postcritical orbit
$$
\widetilde \MM (\rho_0)(x)=\sum_{f(y)=x} \frac{\rho'_{reg}(y) }{|f'(y)|f'(y)}
-\rho_{reg}'(x)\, .
$$
In other words, we have (in $BV_p$)
$$
\widetilde \MM (\rho_0)=\bar \MM(\rho'_{reg})
-\rho_{reg}'\, ,
$$
where
$\bar \MM \varphi(x)= \sum_{f(y)=x} \frac{\varphi(y) }{|f'(y)|f'(y)}$.
So we have reduced the claim (\ref{ccclaim}) to
\begin{align*}
-\alpha \rho_{reg}' 
-  (\id-\LL_0)^{-1}(\id-\Pi_0)(\LL_0-\id)(\alpha \rho'_{reg})
 = 0\, ,
\end{align*} 
that is, using $\Pi_0 \LL_0=\Pi_0$,
$$
(\LL_0 -\id)(\alpha \rho_{reg}') =
(\id-  \Pi_0)(\LL_0-\id)(\alpha \rho'_{reg})=(\LL_0-\id)(\alpha \rho'_{reg}) \, .
$$

Finally, since $\rho_0\in BV_p$,
and since $t\mapsto \nu_t$ and $t\mapsto \tilde \rho_t$
are differentiable in $BV_p$ and $BV_p^*$, respectively, we have 
(in $BV_p$)
\begin{equation}\label{theend}
\partial_t (\nu_t(\rho_0) \tilde \rho_t)|_{t=0} = 
\partial_t (\nu_t(\rho_0))|_{t=0} \rho_0 +
\partial_t ( \tilde \rho_t)|_{t=0}\, .
\end{equation}
Take the Lebesgue average of both sides of (\ref{theend}).
Since $\partial_t \int \tilde \rho_t\, dx=0$
(because each $\tilde \rho_t\, dt$ is a probability),
and since  
$$-\int (\id-\LL_0)^{-1} (X'\rho_{sal}+ (X\rho_{reg})')\, dx=0$$
(use again $\Pi_0 (X'\rho_0+X \rho_{reg}')=0$),
we find that $\partial_t (\nu_t(\rho_0))|_{t=0} \int \rho_0 \, dx=0$.  Therefore, 
$\partial_t (\nu_t(\rho_0))|_{t=0}$, and
putting together (\ref{deccomp}), (\ref{otherterm}),  (\ref{resid}), (\ref{ccclaim}), and
(\ref{theend}),  
we have proved the theorem.
\end{proof}

\smallskip

We have (a simplification of Lemma~\ref{jensenlemma}):
\begin{lemma}\label{gla}
Let $f_t$ be a $C^2$ family of piecewise expanding $C^3$
unimodal maps in the topological class of $f_0$. Set $v=\partial_t f_t|_{t=0}$.
For any $p >1$ the  
map 
$ t \mapsto g_{t}=  \frac{1}{|f'_t \circ h_t|} \in BV_p$ is 
 $C^1$ in a neighbourhood of $0$, and
$\partial_t g_t|_{t=0} =-\frac{f''_0\alpha +v'}{|f_0'|f'_0}$.
\end{lemma}

\begin{proof}
Differentiability  follows from Lemma~\ref{jensenlemma}
applied to $\psi \equiv 0$. The value of the derivative
is given by (\ref{first}--\ref{second}) in the proof of that lemma,
since $\partial_t f'_t|_{t=0}=v'$.
\end{proof}

\begin{remark}
We have the following strengthening of Lemma~\ref{gla}
if $f_t$ is a $C^{3}$ family of piecewise expanding $C^{4}$
unimodal maps:
For any $p >1$ the  
map 
$ t \mapsto g_{t}=  \frac{1}{|f'_t \circ h_t|} \in BV_p$ is 
$C^2$ in a neighbourhood of $0$, and
$| g_t-g_0+t\frac{f''_0\alpha +v'}{|f'_0|f'_0}|=O(t)$.
Recalling  (\ref{deriv0}--\ref{deriv}), this implies that
$ \|\frac{ \tilde \LL_t-\LL_0 }{t} -  \MM  \|_{BV_p} =O(t)$. 
\end{remark}

We shall get the following (new) result as a corollary of Theorems ~\ref{pressuretrick}
and ~\ref{theformula}:

\begin{corollary}\label{thebest}
If $f_t$  is a $C^2$ family of piecewise
expanding $C^3$ unimodal maps
in the topological class of $f_0$, and if
$\partial_t f_t |_{t=0}=X \circ f_0$ for a $C^2$ function $X$,
then  there exists $\epsilon>0$ so that
$t \mapsto \mu_t$ is $C^{1}$ from $(-\epsilon,\epsilon)$ to Radon measures.
\end{corollary}

In particular, under the assumptions of Corollary~\ref{thebest},
the Radon measure 
$$
- \alpha_t  \rho'_{sal,t}- (\id- \LL_t)^{-1} (X_t' \rho_{sal,t}+(X_t\rho_{reg,t})') \, dx
$$
(recall (\ref{formula})) is continuous as a
function of $t$.
(Here, $\alpha_t$ solves (\ref{tceeq}) for $f_t$
and $v_t=\partial_s f_s|_{s=t}$,
and $X_t \circ f_t = v_t$.)
This fact is not clear a priori from the formula.

\begin{remark}We expect that a careful analysis
of the term (\ref{otherterm}) for $C^1$ functions $\psi$ would allow to bypass
the reference to Theorem~\ref{pressuretrick} in the proof of
Corollary~\ref{thebest}.
\end{remark}

\begin{proof}[Proof of Corollary~\ref{thebest}]
We want to show that $t\mapsto \partial_u \mu_u|_{u=t}=\tilde \mu_t$ is continuous: 
We know that $\tilde \mu_t$ exists
for all small $t$
(as a Radon measure) by Theorem~\ref{theformula}.
Clearly, $|\int \psi \, d\tilde \mu_t |\le C \sup |\psi|$ for all continuous $\psi$
and all small enough $t$.

Assume for a contradiction that $t\mapsto \tilde \mu_t$ is discontinuous at  $t_0$.
This means that there exist $\psi\in C^0$, with $\sup |\psi|=1$,
$\delta>0$, and a sequence $t_m$
with $|t_m-t_0|<1/m$, so that 
$|\int \psi \, d\tilde \mu_{t_0}- \int \psi \, d\tilde \mu_{t_m}|>\delta$
for all $m$. 
Take $\tilde \psi\in C^1$ so that $\sup|\psi -\tilde \psi|< \delta/4$.
Then $|\int \tilde \psi \, d\tilde \mu_{t_0}- \int \tilde \psi \, d\tilde \mu_{t_m}|>\delta/2$
for all $m$.
But  Theorem~\ref{pressuretrick} implies
 $|\int \tilde \psi \, d\tilde \mu_{t_0}- \int \tilde \psi \, d\tilde \mu_{t_m}|< \delta$
 if $m$ is large enough, a contradiction.
 \end{proof}

\begin{appendix}
\section{A consequence of the Keller-Liverani bounds from \cite{BS1}}

We state here for the record
an immediate corollary of \cite[Proposition~3.3]{BS1} which was based
on results in \cite{KL}
(see Remark~\ref{famdef} and note that the assumptions
below imply $\sup_I |f_t-\tilde f_t|=O(t^2)$):

\begin{proposition}\label{KLbd}
Let $f_t$ be a $C^2$ family of piecewise expanding $C^2$ unimodal maps.
Assume  that
$f_0$ is mixing and good,
and that $f_t$ is tangent to the topological class of $f_0$,
denoting by $\tilde f_t$ a family in the topological class of $f_0$ with
$\tilde f_0=f_0$ and $\partial f_t|_{t=0}=\partial \tilde f_t|_{t=0}$.

Let $\mu_t=\rho_t \, dx$  and $\tilde \mu_t=\tilde \rho_t \, dx$
be the $SRB$ measures of $f_t$ and $\tilde f_t$, respectively.
Then
for any $\xi <2$ there exists $C>0$  so that for all small $t$
$$
\|\rho_t -\tilde \rho_t\|_{L^1(Leb)} \le C |t|^\xi  \, .
$$
\end{proposition}
\end{appendix}

\bibliographystyle{amsplain}

\end{document}